\documentclass[12pt,a4paper]{article}
\usepackage[cp1251]{inputenc}
\usepackage[russian]{babel}
\usepackage{indentfirst}
\usepackage{amsmath}
\usepackage{graphicx}
\usepackage{float}
\usepackage{stmaryrd}
\usepackage{bbm}
\usepackage[all]{xy}
\usepackage{amsfonts}
\usepackage{amssymb}

\usepackage{amscd}
\usepackage{mathrsfs}
\usepackage{tipa}
\usepackage{mathtext,amsmath, amstext, amsgen, amsthm}
\usepackage{amsopn, amsfonts, euscript,amscd,amssymb,latexsym}
 \newtheorem{theorem}{Theorem}

 \newtheorem{prop}{Proposition}
 \theoremstyle{definition}
 \newtheorem{defn}{Definition}
 \theoremstyle{definition}

 \numberwithin{equation}{section}

\begin{document}

\begin{center}
{\bf \large BASIC AUTOMORPHISM GROUPS OF CARTAN  \\

\vspace{2mm}
 FOLIATIONS COVERED BY FIBRATIONS}\\
\end{center}

\begin{center}
{K.I. SHEINA}\\
\end{center}

\begin{center}
{\it National Research University Higher School of Economics, \\
Department of Informatics, Mathematics and Computer Science,\\
ul. Bolshaja Pecherskaja, 25/12, Nizhny Novgorod, 603155, Russia\\
} 
\,\,\,E-mail:{ksheina@hse.ru} 
\end{center}

\vspace{4mm} {\bf Abstract} The basic automorphism group ${A}_B(M,F)$ of a Cartan foliation $(M, F)$ is the quotient group of the  automorphism group of $(M, F)$ by the normal subgroup, which preserves every leaf invariant. For Cartan foliations covered by fibrations, we find sufficient conditions for the existence of a structure of a finite-dimensional Lie group in  their   basic automorphism groups. Estimates of the dimension of these groups are obtained. 
For some class of Cartan foliations with integrable an Ehresmann connection, a method for finding  of basic automorphism groups is specified.

{\it \bf Keynwords: }{foliation; Cartan foliation; Lie group; basic automorphism; auto\-mor\-phism group}


\section{Introduction}

 One of the main objects associated with a geometric structure on a smooth manifold is its automorphism group. In the introduction to the monograph by S. Kobayashi~\cite{K}, it was emphasized that the existence of a structure of a finite-dimensional Lie group in the group of automorphisms of a manifold with a geometric structure is one of the central problems in differential geometry.
 
 As is known, the solved Hilbert's 5-th problem is devoted to finding conditions under which a topological group admits the structure of a Lie group \cite{Skl}. It is known from the numerous works of E. Cartan, R. Mayer, H. Steenrod, K. Nomizu, S. Kabayashi, S. Ehresmann and other authors that the automorphism groups of many geometric structures are Lie groups of transformations (see 
 overview \cite{Chu}).
 
 The spaces that are now called Cartan geometries were introduced by E. Cartan in the 1920-s. The theory of Cartan geometries is presented in the monographs of A.~$\check{C}$ap, J.~Slovak \cite{C-S}, R.V. Sharpe \cite{Shar}, M. Krampin and D. Saunders \cite{Cr}. Currently, Cartan geometries and Cartan foliations are studied by many mathematicians and find application in various physical theories, see, for example, \cite{Baz}, \cite{CGH}, \cite{Pec} and \cite{J}, \cite{Z}.
 
 Let $(M, F)$ be a smooth foliation. Recall that  geometry structure on the  manifold $M$ is called   transverse to  $(M,F)$ if it is  a invariant with respect to local holonomic diffeomorphisms. Another, equivalent definition of a transverse geometric structure, which is represented by Cartan geometry, is given in Section \ref{S2}. Morphisms are under\-stood as local diffeomorphisms mapping leaves onto leaves and preserving transverse geometries (the precise definition see in Section \ref{S2}).  Let us denote by $\mathfrak C\mathfrak F$ the category of Cartan foliations.

This paper is devoted to the investigation of automorphism groups of Cartan foliations, i.e. foliations that admit Cartan geometries as  transverse structures. The study of Cartan foliations is motivated by the fact that such broad classes of foliations as parabolic, conformal, projective, pseudo-Riemannian, Lorentzian, Weyl, transverse homogeneous foliations and foliations with transverse linear connection belong to Cartan foliations. Therefore, the investigation  of Cartan foliations allows us to study the general properties of these foliations from a single point of view, while many authors study them separately.
 
Let us denote by $A(M, F)$ the  group of all the automorphisms of the
Cartan foliation $(M, F)$ in the category $\mathfrak{C}\mathfrak{F}$. The group
\begin{equation}
	{ A}_L(M,F):=\{f\in{A}(M, F)\mid
	f(L_\alpha)=L_\alpha\,\,\,\forall L_\alpha\in F\} \nonumber
\end{equation}
is a normal subgroup of the group ${A}(M, F)$ and called {\it the group of leaf automorphisms} of   $(M, F).$
The quotient group ${ A}(M,F)/{A}_L(M,F)$ is called  {\it the basic automorphism group } and  denoted by ${A}_B(M,F).$

We study the groups of basic automorphisms ${A}_B(M,F)$ of Cartan foliations $(M,F)$ covered by fibration and find sufficient conditions for the existence of a structure of a finite-dimensional Lie group in the group ${A}_B(M,F)$. J. Leslie \cite{Les} was the first who  solved a similar problem for smooth
foliations on compact manifolds and considered an application to foliations with  transverse G-structures. For foliations with complete
transversely projectable affine connec\-tion, this problem was raised
by I.V. Belko~\cite{Bel}. Foliations $(M,F)$ with effective transverse rigid geometries were
investigated by\- N.I.~Zhu\-kova \cite{ZhR} where an algebraic
invariant $\mathfrak g_0 = \mathfrak g_0(M, F),$ called the structural Lie
algebra of $(M,F)$, was constructed and it was proved that $\mathfrak g_0 = 0$
is a sufficient condition for the  existence of a unique  Lie
group structure in the basic automorphism group of this
foliation.   In \cite{SZ}, the existence of a Lie group structure was investigated  in the  basic automorphism groups of Cartan foliations modeled on inefficient Cartan geometries.

\section{Main results}

Among the Cartan foliations, foliations covered by fibrations are distinguished. 

\begin{defn}\label{opr1} Let $\kappa: \widetilde{M}\to M$ be the universal covering map.
We say that a smooth foliation $(M, F)$
{\it is covered by fibration} if the induced foliation $(\widetilde{M}, \widetilde{F})$
is formed by fibres of a locally trivial fibration
$\widetilde{r}: \widetilde{M}\to B.$
\end{defn}
The following theorem describes the global structure of  Cartan foliations covered by fibrations.
\begin{theorem}\label{Th2}
Let $(M, F)$ be a  Cartan foliation modeled on a Cartan geometry $\xi$ covered by the fibration
$\widetilde{r}:\widetilde{M}\to B$, where $\widetilde{\kappa}:
\widetilde{M}\to M$  is the universal covering map. Then:
\begin{enumerate} \itemsep=0pt
\item[(1)] there exists a regular covering map $\kappa: \widehat{M}\to M$
 such that the induced foliation $\widehat{F}$ is made up of fibres of the locally
trivial bundle $r: \widehat{M}\to B$ over a simply connected 
manifold $B$, and $\xi$ induces on $B$ a Cartan geometry $\eta$ that is locally isomorphic to $\xi$;
\item[(2)] an epimorphism $\chi: \pi_1(M, x)\to \Psi,\, x\in M,$
of the fundamental group $\pi_1(M, x)$ onto a subgroup $\Psi$ of the automorphism Lie group  $Aut(B,\eta)$ of the Cartan manifold $(B,\eta)$ is determined;
\item[(3)] the group of deck transformations of  the covering $\kappa: \widehat{M}\to M$ is isomorphic to the group $\Psi$.
\end{enumerate}
\end{theorem}
\begin{defn} The group $\Psi = \Psi(M, F)$ satisfying Theorem~\ref{Th2}
is called the {\it global holonomy group} of the Cartan foliation $(M, F)$
covered by fibration.
\end{defn}

We give a detailed proof of the following theorem, formulated without a proof in the work \cite[Prop.~8]{SZ}. Theorem \ref{Th5} establishes a connection between the basic automorphism group  $A_{B}(M,F)$ of a Cartan foliation $(M,F)$ covered by fibration and its global holonomy group $\Psi$.

\begin{theorem}\label{Th5}
Let $(M, F)$ be a Cartan foliation  covered by  fibration
$r:\widehat{M}\rightarrow B$, and $B$ is the simply connected
Cartan manifold. Suppose that the
global holonomy group $\Psi$ is a discrete subgroup of the Lie group
$Aut(B,\eta)$. Let $N(\Psi)$ be the normalizer of $\Psi$ in $Aut(B,
\eta)$. Then the basic automorphism group $A_B(M, F)$ is a Lie
group which  is isomorphic to an open-closed subgroup of the Lie
quotient group $N(\Psi)/\Psi$ and $\dim(A_{B}(M,
F))=\dim(N(\Psi)/\Psi).$  The structure of the Lie group in  $A_{B}(M,F)$ is unique.
\end{theorem}

The following theorem specifies a method for computing   basic automorphism groups for Cartan foliations with an integrable Ehresmann connection.

\begin{theorem}

\label{Th7}

 Let $(M,F)$ be Cartan foliation with an integrable Ehresmann connection. Then

1.There is a regular cover $\kappa: \widetilde{M}\to M$ such that $\widetilde{M}~=L_0~\times~B,$  where $L_0$ is a manifold diffeomorphic to any leaf with a trivial holonomy group and B is a simply connected manifold, and the induced foliation $\widetilde{F}=\kappa^*F$ is formed by   leaves of the canonical projection $r: L_0\times B\to B$ onto the second factor, and  Cartan geometry $\eta$   is induced on $B$, with respect to which $\kappa$  is a morphism of Cartan foliations $(M, F)$ and $(\widetilde{M}, \widetilde{F})$ in the category $\mathfrak
C\mathfrak F$.

2. The foliation  $(M,F)$ is an $(Aut(B,\eta),B)$-foliation.

3. If moreover, the normalizer $N(\Psi)$ of global holonomy group $\Psi$ is
equal to  the centralizer $Z(\Psi)$ of $\Psi$ in the group $Aut(B,\eta)$, then
$$A_B(M, F)\cong N(\Psi)/\Psi.$$
\end{theorem}

Using Theorem \ref{Th7}, we construct an example of computing the  basic auto\-morphism group of some conformal foliation of an arbitrary codimension $q$, where $q\geq 3$, on a $(q+1)$-dimen\-sional manifold in Section \ref{Ex}. Some other examples are constructed in~\cite{SZ}.

\section{The category of Cartan foliations}
\label{S2}
\paragraph{The category of Cartan geometries}

Let $G$ and $H$  be  Lie groups with the Lie algebras  $\mathfrak{g}$ and $\mathfrak{h}$ relatively.  Let $H$ be a closed subgroup of $G$. {\it A Cartan geometry} of type  $(G,H)$ on the smooth manifold $N$  is a principal
$H$-bundle $P(N,H)$ with
a $\mathfrak{g}$-valued $1$-form $\omega$ on $P$ satisfying the following  conditions:\\
($c_1$) the map $\omega_{u}:T_{u}P\rightarrow \mathfrak{g}$ is an isomorphism of
vector spaces for every $u\in P$;\\
($c_{2}$)  $\omega({A^{*}})=A$ for every $A\in\mathfrak{h}$, where $A^{*}$
is the fundamental vector field determined by $A$;\\
($c_{3}$) $R^{*}_{h}\omega=Ad_{G}(h^{-1})\omega$ for every $h\in H$,
where $Ad_{G}:H\rightarrow GL(\mathfrak{g})$  is the adjoint representation
of the Lie subgroup $H$ of $G$ in the Lie algebra $\mathfrak{g}$.

 The $\mathfrak{g}$-valued form $\omega$ is called a {\it Cartan
 connection form}. This Cartan geometry is denoted  by $\xi=(P(N,H),\omega)$.
 The pair $(N,\xi)$ is called a {\it Cartan manifold}.

Maximal normal subgroup $K$ of the group
$G$ belonging to $H$ is called the {\it kernel} of pair $(G,H)$.
We denote the Lie algebra  of the  group  $K$  by $\mathfrak{k}$. The Cartan geometry $\xi=(P(M,H),\omega)$ of 
 type $(G,H)$ is called {\it effective}
if the kernel $K$ of the pair $(G,H)$ is trivial. Further, we assume that all Cartan
 geometries under consideration are effective.

 Let $\xi=(P(N,H),\omega)$ and  $\xi'=(P'(N',H),\omega')$ be two
Cartan geometries with the same structure Lie group $H$.
The smooth map $\Gamma:P\to P'$ is called a morphism from $\xi$ to $\xi'$ if
$\Gamma^{*}\omega'=\omega$ and $R_{a}\circ\Gamma=\Gamma\circ R_{a}\,\,\, \forall a\in H$.
 The category of Cartan geometries is denoted  by ${\mathfrak
	C }{\mathfrak a}{\mathfrak r}$. If $\Gamma\in Mor(\xi, \xi')$, then the projection 
$\gamma:N\to N'$ is defined such that $p'\circ \Gamma=\gamma\circ p,$  where $p:P\to N$ and $p':P'\to N'$
 are the projections of the respective $H$-bundles.

The projection $\gamma$ is called {\it an automorphism of the Cartan
manifold $(N,\xi)$}. Denote by $Aut(N,\xi)$ the full automorphism
group of a Cartan foliation $(N,\xi)$ and by $Aut(\xi)$ the full automorphism group of a Cartan geometry
$\xi$. Let $A(P,\omega):=\{\Gamma\in Diff
(P)\,|\\ \,{\Gamma^{*}\omega=\omega}\}$ be the automorphism group  of
the parallelizable manifold $(P,\omega)$, which is known to be a Lie group, and $\dim(A(P,\omega)) \leq \dim P.$

Remark, that $Aut(\xi) =\{\Gamma\in A(P,\omega)\,| \,
\Gamma\circ R_{a}=R_{a}\circ \Gamma\, \, \forall a\in H\}$  is a closed Lie subgroup of the Lie group $A(P,\omega)$. 
Therefore, $Aut(\xi)$ is a  Lie group, and due to the effictivity of a Cartan geometry $\xi$, there exists a Lie group isomorphism $$\sigma:A^{H}(P,\omega)\to Aut(N,\xi):\Gamma\mapsto \gamma$$
mapping $\Gamma\in A^{H}(P,\omega)$ to its projection $\gamma$.

\paragraph {Cartan foliations}

Let $N$ be a smooth $q$-dimensional manifold, the connectivity of
which is not assumed. Let $M$ be a smooth $n$-dimensional manifold, where $0<q<n$. 
Assume, that $\xi =(P(N, H),\omega)$ is a Cartan geometry of type $(G,H)$  on the manifold $N$. Let  $p:P\to N$ be the projection the principal $H$-bundle.
For every open subset $V\subset N$, the  Cartan geometry
$\xi_{V}=(P_{V}(V, H), \omega_{V})$ of the same type $(G,H)$ is induced, where 
$P_{V}:=p^{-1}(V)$ and $\omega_{V}:=\omega|_{P_V}$.
Remind that  $(N, \xi)$-\textit{cocycle} on $M$ is a
family $\zeta=\{U_{i},f_{i},\{\gamma_{ij}\}\}_{ij\in J}$ satisfying the following conditions:\\
1)the set  $\{U_{i}\,|\, i\in J\}$ is a covering of the manifold $M$ by open connected  subsets $U_{i}$ of $M$, and every $f_{i}: U_{i}\to N$
is a submersion with connected fibres;\\
2) if $U_{i}\cap U_{j}\neq\emptyset,\, i,j\in J$, then there exists an isomorphism $\Gamma_{ij}:\xi_{f_{j}(U_{i}\cap U_{j})}\to \xi_{f_{i}(U_{i}\cap U_{j})}$
of the  Cartan geometries induced on open subsets ${f_{j}(U_{i}\cap U_{j})}$ and ${f_{i}(U_{i}\cap U_{j})}$ such that
the projection $\gamma_{ij}$ of the isomorphism  $\Gamma_{ij}$ satisfies the equality $f_{i}=\gamma_{ij} \circ f_{j}$ on $U_{i}\cap U_{j},\, i,j\in J$;\\
3) if $U_{i}\cap U_{j}\cap U_{k}\neq\emptyset$, then $\gamma_{ij}\circ \gamma_{jk}=\gamma_{ik}$ 
for all $x\in U_i\cap U_j\cap U_k$ and $\gamma_{ii}=id_{U_{i}}$, $i, j, k\in J$.

Two $N$-cocycles are called  {\it equivalent} if there exists  an
$N$-cocycle containing both of these cocycles. Let $[\{U_{i},f_{i},\{\gamma_{ij}\}\}_{ij\in J}]$ be the
equivalence class of $N$-cocycles on the  manifold $M$ containing the
cocycle  $\zeta=\{U_{i},f_{i},\{\gamma_{ij}\}\}_{ij\in J}$. Denote by $\Sigma$
the set of fibres  of all the submersions $f_{i}$ of this equivalence
class. Note, that $\Sigma$ is the base of some new topology $\tau$
on $M$. The path-connected components of the topological space $(M,
\tau)$ form a partition $F:=\{L_{\alpha}\,|\,\,\alpha\in
\mathfrak{J}\}$ of the manifold $M$. The pair $(M,F)$ is called {\it a
 Cartan foliation of codimension $q$ } modeled on the Cartan geometry
  $\xi=(P(N,H),\omega)$ which is called {\it a transverse Cartan geometry} for 
   $(M,F)$. Subsets $L_{\alpha}, \, \alpha\in \mathfrak{J}$ are called  {\it leaves} of this foliation. It is said that $(M,F)$ is given by the  $(N,\xi)$-cocycle $\zeta$.

\paragraph{Morphisms in the category of Cartan foliations}
Let $(M,F)$ and $(M',F')$ are Cartan foliations defined by an
$(N,\xi)$-cocycle $\zeta=\{U_i,f_i, \{\gamma _{ij}\}\}_{ij\in J}$ and an
$(N',\xi')$-cocycle $\zeta'=\{U'_r,f'_r, \{\gamma'_{rs}\}\}_{i'j'\in J'}$,
respectively. All objects be\-long\-ing to $\zeta'$ are dis\-tingui\-shed by prime.
Let $f\colon M\to M'$ be a smooth map which is a local
isomorphism in the foliation category ${\mathfrak F\mathfrak o\mathfrak l}.$
Hence for any $x\in M$ and $y:=f(x)$ there exist neighborhoods
$U_k\ni x$ and $U'_s\ni y$ from $\zeta$ and $\zeta'$ respectively,
 a diffeomorphism $\varphi\colon V_k\to V'_s,$ where
$V_k:=f_k(U_k)$ and $V'_s:=f'_s(U'_s),$ satisfying the relations
$f(U_k)=U'_s$ and $\varphi\circ f_k=f'_s\circ f|_{U_k}$. Further we shall
use the following notations: $P_{k}:=P|_{V_{k}},\,\,\,P'_{s}:=P'|_{V'_{s}}$
and $p_{k}:=p|_{P_{k}},\,\,$ $p'_{s}:=p|_{P'_{s}}$

We say that $f$ preserves  transverse Cartan geometry if every such
dif\-feo\-mor\-phism $\varphi\colon V_k\to V'_s$ is an isomorphism of
the induced Cartan geometries $(V_k,\xi_{V_k})$ and
$(V'_s,\xi'_{V'_s})$. This means the existence of isomorphism
$\Phi: P_{k}\to P'_{s}$ in the category $\mathfrak C\mathfrak a\mathfrak r$
with the projection $\varphi$ such that  the following diagram
\begin{equation}
\xymatrix {&{P_{k}}\ar[d]_{p_{k}}\ar[rd]^{\Phi}&\\
 M\supset{U_k}\ar[r]^{f_{k}}\ar[rd]^{f|_{U_{k}}} & V_{k}\ar[rd]^{\varphi}  &P'_{s}\ar[d]^{p'_{s}}\\
& M'\supset{U'_{s}}\ar[r]^{f'_{s}}& V'_{s}}\nonumber
\end{equation}
 is commutative. We emphasize that the indicated above isomorphism
$\Phi:P_{k}\rightarrow P'_{s}$  is  unique if the transverse Cartan
geometries  are  effective. The introduced  concept  is well
defined, i.~e., it does not depend of the choice of neighborhoods
$U_k$ and $U'_k$ from the cocycles $\zeta$ and $\zeta'.$

\begin{defn}
By a {\it morphism of two Cartan foliations $(M, F)$ and $(M', F')$
} we mean a local diffeomorphism $f:M\to M'$  which transforms leaves to leaves and preserves
transverse Cartan structure. The category $\mathfrak C\mathfrak{F}$
objects of which are Cartan foliations, mor\-phisms are their
mor\-phisms, is called {\it the category of Cartan foliations.}
\end{defn}

\section{Ehresmann connections for foliations}\label{ssEr}
R. A. Blumenthal and J. J. Hebda \cite{BH} introduced the notion of Ehresmann connection for foliation $(M,F)$ as a natural generalization of Ehresmann connection for submersions.

Let $(M,F)$ be a foliation of codimension $q$ and $\mathfrak M$ be a
smooth $q$-dimensional distribution on $M$ that is transverse to
the foliation $F,$ i. e. $T_{x}M=\mathfrak{M}_{x}\oplus T_{x}F\,\,\,\forall x\in M$.  The piecewise smooth integral curves of the
distribution $\mathfrak M$ are said to be {\it horizontal,} and the
piecewise smooth curves in the leaves are said to be {\it vertical}. A
piecewise smooth mapping $H$ of the square $I_1\times I_2$ to $M$
is called a {\it vertical-horizontal homotopy} if the curve
$H|_{\{s\}\times I_2}$ is vertical for any fixed $s\in I_1$ and the curve
$H|_{I_1\times\{t\}}$ is horizontal for any fixed $t\in I_2.$ In this
case, the pair of paths $(H|_{I_1\times \{0\}},H|_{\{0\}\times
	I_2})$  is called the {\it base} of $H.$ It is well known that there
exists at most one vertical-horizontal homotopy with a given base.

A distribution $\mathfrak M$ is called an {\it Ehresmann
	connection for a foliation} $(M,F)$ (in the sense of R. A. Blumenthal
and J. J. Hebda \cite{BH}) if, for any pair of paths $(\sigma, h)$ in $M$ with
a common initial point $\sigma(0) = h(0),$ where $\sigma$ is a
horizontal curve and $h$ is a vertical curve, there exists a
vertical-horizontal homotopy $H$ with the base $(\sigma, h).$

\paragraph{A simple foliation with an Ehresmann connection} Let $f:M\to N$ be a submersion with connected fibers. Recall that the foliation $F=\{p^{-1}(z),\, z\in N\}$ formed by the fibers of the submersion is called {\it a simple foliation}. Let $(M,F)$ be an arbitrary smooth foliation with  the Ehresmann connection. It easy to show that, the existens of a covering $\widehat{k}:\widehat{ M}\to M$ such that the lifted foliation is simple implies that the foliation $(M,F)$ is covered by fibration.

\section{Classes of foliations covered by fibrations}
\paragraph{ $(G,X)$-foliations with an Ehresmann connection} Let $X$ be a smooth connected manifold and $G$ be the Lie group of diffeomorphisms of $X$. Recall that the action of a group $G$ on a manifold $X$ is called {\it quasi-analytically} if for any open subset $U\subset X$ and an
element $g\in G$ the equality $g|_U = id_U$ implies
$g = e$, where $e=id_X$.

Assume that a Lie group  $G$ of diffeomorphisms of
	a manifold $X$ acts on $N$ quasi-analytically. A foliation $(M, F)$
	defined by an $X$-cocycle $\{U_i,f_i,\{\gamma_{ij}\}\}_{i,j\in J}$
	is called a $(G,X)$-foliation if for any $U_i\cap U_j
	\neq\emptyset$, $i,j\in J$, there is an element $g\in G$
	such that $\gamma_{ij} = g|_{f_j(U_i\cap U_j)}$.  If, moreover, $(X,\xi)$ is a Cartan manifold and the group $G$ is a subgroup of the automorphism Lie group  $Aut(X,\xi)$, then $(M,F)$ is a Cartan $(G,X)$-foliation.  It follows from [7, Section VI] that every Cartan $(G,X)$-foliations  with Ehresmann connections is a foliation covered by fibration. 
	
\paragraph{Cartan foliation  with a vanishing transverse curvature}
Let $(M,F)$ be a Cartan foliation of type $(G,H)$ with an Ehresmann connection. As is known \cite[Section VI]{Baz}, if the transverse curvature of $(M, F)$
vanishes, then foliation $(M, F)$ is covered by fibration. Consequently, all the  obtained results are valid for Cartan foliations  with zero transverse curvature that admiting an Ehresmann connection.
\paragraph{Conformal foliations of codimension $q, q \geq 3$}
 According to \cite[Thm.~5]{ZhG}, any  non-Riemannian conformal foliation of codimension $q \geq 3$ with an Ehresmann connection  is covered by fibration. 

 \paragraph{Foliations with an integrable Ehresmann connection} Recall that an Ehresmann connection $\mathfrak M$ for a foliation $(M,F)$ is called integrable if the  distribution $\mathfrak M$  is integrable i.e. if there exists the foliation such that $TF^{t}=\mathfrak M$.  Accoding to Kashiwabara's theorem \cite{Kas}, foliations with an integrable  Ehresmann connection are covered by fibrations.

\paragraph{Suspended foliations}
The construction of a suspension  foliation was proposed by A. Haefliger and described in detail in \cite{Min}. Note that suspension  foliations form a class of foliations with integrable Ehresmann connection and are covered by  fibrations.

 \paragraph{Cartan foliation of codimension $q=1$}
Any smooth one-dimensional distribution is integrable, so a Cartan foliation $(M,F)$ of codimension $q=1$  with an Ehresmann connection is covered by fibration.

\section{Proof of Theorems \ref{Th2} and \ref{Th5}}

\subsection{ Regular covering maps}
\begin{defn}
	Let $f:M\rightarrow B$ be a submersion. It is said that $\widehat{h}\in Diff (M)$
	lying over $h\in Diff (B)$ relatively $f$ if $h\circ f=f\circ \widehat{h}$. In this case $\widehat{h}$ is called {\it a lift} of $h$ with respect to $f:M\rightarrow B$.
\end{defn}

Let $\widetilde{\kappa}:\widetilde{{K}}\to K$  be the universal covering map, where $K$ and $\widetilde{{K}}$
are smooth manifolds. By analogy with Theorem~28.10 in \cite{Bus}, it is easy to show
that for any $h~\in~Diff(K)$ there exists $\widetilde{h}~\in~ Diff(\widetilde{K})$
lying over $h$. For an arbitrary covering map the same statement is  incorrect, in general. It is not difficult to prove the following criterion for the existence of lifts of arbitrary diffeomorphisms with respect to regular covers.

\begin{prop}\label{pr8.2}
	Let ${\kappa}:\widehat{K}\to K$ be a smooth regular  covering map with
	the deck transformation group $\Gamma$. A diffeomorphism $\widehat{h}\in Diff(\widehat{K})$	lies over some diffeomorphism  $h\in Diff(K)$ if and only if it satisfies the equality $\widehat{h}\circ\Gamma=\Gamma\circ \widehat{h}$.
\end{prop}

\subsection{Proof of Theorem \ref{Th2}} Suppose that a Cartan foliation $(M, F)$ modeled on an effective
Cartan geometry $\xi = (P(N,H),\omega)$ is covered by a fibration
$\widetilde{r}: \widetilde{M}\to B$, where $\widetilde{\kappa}:
\widetilde{M}\to M$ is the universal covering map. The fibration
$\widetilde{r}: \widetilde{M}\to B$ has connected fibres and simply
connected space $\widetilde{M}$. Therefore, due to the application of
the exact homotopic sequence for this fibration we obtain that the base
manifold $B$ is also simply connected.

For an arbitrary point $b\in B$ take $y\in\widetilde{r}^{-1}(b)$
and $x=\widetilde{\kappa}(y)$. Without loss genera\-lity, we assume
that there is a neighbourhood $U_i$, $x\in U_i$, from the $(N,\xi)$-cocycle\-
$\{U_i,f_i,\{\gamma_{ij}\}\}_{i,j\in J}$ which defines $(M, F)$ and
a neighbourhood $\widetilde{U}_i$, $y\in\widetilde{U}_i$, such that
$\widetilde{\kappa}|_{\widetilde{U}_i}: \widetilde{U}_i\to U_i$ is a
diffeomorphism.

Let $\widetilde{V}_i: = \widetilde{r}(\widetilde{U}_i)$. Then there
exists a diffeomorphism $\phi: \widetilde{V}_i\to V_i$  satisfying
the equality $\phi\circ\widetilde{r}|_{\widetilde{U}_i} =
f_i\circ\widetilde{\kappa}|_{\widetilde{U}_i}$. The diffeomorphism
$\phi$ induces the Cartan geometry $\eta_{\widetilde{V}_i}$ on
$\widetilde{V}_i$ such that $\phi$ becomes the isomorphism
$(\widetilde{V}_i,\eta_{\widetilde{V}_i})$ and $(V_i,\xi_{V_i})$ in
the category $\mathfrak C\mathfrak a \mathfrak r$ of Cartan geometries. The
direct check shows that by this way we define the Cartan geometry
$\eta$ on $B$, and $\eta|_{\widetilde{V}_i} =
\eta_{\widetilde{V}_i}$, $i\in J$. Thus, the statement $(1)$ is proved.

Let us fix points $x_0\in M$ and
$y_0\in\widetilde{\kappa}^{-1}(x_0)\in\widetilde{M}$. Then the fundamental group
$\pi_1(M,x_0)$ acts on the universal covering space $\widetilde{M}$
as the deck transformation group $\widetilde{G}\cong\pi_1(M,x_0)$ of
$\widetilde{\kappa}$. Since $\widetilde{G}$ preserves the inducted
foliation $(\widetilde{M},\widetilde{F})$ formed by fibres of the
fibration $\widetilde{r}: \widetilde{M}\to B$, then every
$\widetilde{\psi}\in\widetilde{G}$ defines $\psi\in Diff(B)$
satisfying the relation $\widetilde{r}\circ\widetilde{\psi} =
\psi\circ\widetilde{r}$. The map $\chi: \widetilde{G}\to\Psi:
\widetilde{\psi}\to\psi$ is a group epimorphism and the statement (2) is proved.

 Observe that $\widetilde{G}$
is a subgroup of the automorphism group $Aut(\widetilde{M}, \widetilde{F})$ of
$(\widetilde{M}, \widetilde{F})$ in the category $\mathfrak C\mathfrak F$.
Therefore $\Psi$ is a subgroup of the automorphism group
$Aut(B,\eta)$ in the category  of Cartan
geometries $\mathfrak C\mathfrak a\mathfrak r$. The kernel $ker(\chi)$ of $\chi$
determines the quotient manifold $\widehat{M}: =
\widetilde{M}/ker(\chi)$ with the quotient map
$\widehat{\kappa}: \widetilde{M}\to\widehat{M}$ and the quotient
group $\widehat{G}: = \widetilde{G}/ker(\chi)$ such that $M
\cong\widehat{M}/\widehat{G}$. The quotient map $\kappa: \widehat{M}\to M$ is
the required regular covering map, with $\widehat{G}$ acts on $\widehat{M}$ as
a deck transformation group. The map $\widehat{G}\to\Psi:
\widetilde{\psi}\cdot ker(\chi)\mapsto\chi(\widetilde{\psi})$, $\widetilde{\psi}\in\widetilde{G}$,
is a group isomorphism. Thus the statment  (3) is proved.

\subsection{The associated foliated bundle}
	Let $(M, F)$ be a Cartan foliation modeled on Cartan geometry
 $\xi=(P({\mathcal {N}},H),\omega)$ of type  $(G,H)$.
	Then there exists a principal $H$-bundle with the projection
	$\pi :\mathcal{R}\to M$, the $H$-invariant foliation $(\mathcal{R},\mathcal {F})$ and
the	$\mathfrak g$-valued $H$-equivariant $1$-form $\beta$ on ${\mathcal{R}}$ which
	satisfy the following conditions:
	\begin{enumerate}
		\item [(i)] $\beta(A^*)=A$ for any $A\in\mathfrak h$;
		\item [(ii)] the mapping $\beta_u:T_u {\mathcal{R}}\to \mathfrak g$ $\forall u\in {\mathcal{R}} $
		is surjective, and  $ker(\beta_u) = T_u{\mathcal{F} }$;
		\item [(iii)] the foliation $({\mathcal{R}},\mathcal{F})$ is transversely parallelizable;
		\item [(iv)]  the Lie derivative $L_X{\beta}$ is equal to zero for every vector field
		$X$ tangent to the foliation $({\mathcal{R}}, \mathcal{F}).$
	\end{enumerate}

\begin{defn} The principal $H$-bundle ${\mathcal{R}}(M,H)$ 
	is called  {\it the foliated bundle} over the Cartan foliation $(M,F)$. The foliation
	$({\mathcal{R}}, \mathcal{F})$ is called {\it the lifted foliation} for the Cartan foliation $(M,F).$
\end{defn}

If the lifted foliation $(\mathcal {R}, \mathcal {F})$ is formed by fibres
of the locally trivial fibration $\pi_b:~{\mathcal {R}}~\to~W$, then $W=\mathcal {R}/\mathcal {F}$ is a smooth manifold , and a $\mathfrak{g}$-valued\- $1$-form $\beta$ such that $\pi_b^*\widetilde{\beta}: = \beta$ and  locally free action of the Lie group $H$ on $W$ are induced. In this case, $(W,\widetilde{\beta})$ is a parallelizable manifold and $A(W,\widetilde{\beta})$ is the Lie group of its automorphisms that acts freely on $W$. Further, as above, by ${A}^{H}(W,\widetilde{\beta})$ we denote the closed Lie subgroup of $A(W,\widetilde{\beta})$ formed by transformations commuting with the induced action of the Lie group $H$ on $W$.
\qed

\subsection{Proof of Theorem \ref{Th5}}
Suppose that a Cartan foliation $(M, F)$ is covered by fibration. By definition $\ref{opr1}$, the induced  foliation  $(\widetilde{M}, \widetilde{F})$
on the space of the universal  covering $\widetilde{\kappa}:\widetilde{M}\to M$
is defined by  a locally trivial fibration $\widetilde{r}:\widetilde{M}\to B$.
Due to Theorem~\ref{Th2}, the regular covering map $\kappa:\widehat{M}\rightarrow M$
and locally trivial fibration $r: \widehat{M} \rightarrow B$ are defined, where $B$
is a simply connected manifold with the inducted Cartan geometry $\eta$. Let $\Psi$
be the global holonomy group of the foliation $(M,F)$, then $\Psi$ is isomorphic to the
deck transformations group $G$ of the covering $\kappa:\widehat{M}\rightarrow M$. Since the
manifold $\widetilde{M}$ is simply  connected, there exists the universal
covering map $\widehat{\kappa}:\widetilde{M}\to \widehat{ M}$ satisfying the equality
$\kappa\circ \widehat{\kappa}=\widetilde{\kappa}$.  Let $\widetilde{G}$, $G$ and $\widehat{G}$  be the deck transformation groups  of the covering maps
$\widetilde{\kappa}$,  $\kappa$ and $\widehat{\kappa}$ respectively,
with $\Psi\cong G\cong\widetilde{G}/\widehat{G}$.

Let  us consider  the following preimages of the $H$-bundle ${\mathcal {R}}$ respectively
$\widetilde{\kappa}$ and  $\kappa$
\begin{equation}
	\widetilde{\mathcal {R}}:=\{(\widetilde{x},u)\in \widetilde{M}\times {\mathcal {R}}\,|\,\widetilde{\kappa}(\widetilde{x})=\pi(u)\}=\widetilde{\kappa}^{*}{\mathcal {R}}\,\,\, and \nonumber
\end{equation}
\begin{equation}
	\widehat{\mathcal {R}}:=\{(\widehat{x},u)\in \widehat{M}\times {\mathcal {R}}\,|
	\,\kappa(\widehat{x})=\pi(u)\}=\kappa^{*}{\mathcal {R}}. \nonumber
\end{equation}

Remark that the maps
\begin{equation}
	\widetilde{\theta}:\widetilde{{\mathcal {R}}}\to {\mathcal {R}}:(\widetilde{x}, u)\mapsto (\widetilde{\kappa}(\widetilde{x}), u)\,\,\,\,\, \forall(\widetilde{x}, u)\in \widetilde{{\mathcal {R}}},\nonumber
\end{equation}
\begin{equation}
	\theta:\widehat{{\mathcal {R}}} \to {\mathcal {R}}:(\widehat{x}, u)\mapsto (\kappa(\widehat{x}),u) \,\,\,\,\, \forall(\widehat{x}, u)\in \widehat{{\mathcal {R}}},\nonumber\\
\end{equation}
\begin{equation}
	\widehat{\theta}:\widetilde{{\mathcal {R}}}\to \widehat{{\mathcal {R}}}:(\widetilde{x}, u)\mapsto (\widehat{\kappa}(\widetilde{x}), u)\,\,\,\,\, \forall(\widetilde{x}, u)\in \widetilde{{\mathcal {R}}},\nonumber
\end{equation}
are regular covering maps  with the deck transformation groups
$\widetilde{\Gamma}$, $\Gamma$  and $\widehat{\Gamma}$, respectively,
which are isomorphic to the relevant groups $\widetilde{G}$, $G$ and
$\widehat{G}$, i.e. $\widetilde{\Gamma}\cong \widetilde{G}$,
$\Gamma\cong G$ and $\widehat{\Gamma}\cong \widehat{G}$.

Let $(\widetilde{\mathcal {R}}, \widetilde{\mathcal {F}})$ and
$(\widehat{\mathcal {R}},\widehat{\mathcal {F}})$ be the corresponding lifted
foliations for $(M,F)$. Since $(\widetilde{M}, \widetilde{F})$ and $(\widehat{M}, \widehat{F})$
are simple foliations, then $(\widetilde{\mathcal {R}},\widetilde{\mathcal {F}})$ and
$(\widehat{\mathcal {R}},\widehat{\mathcal {F}})$ are also simple foliations,
which are formed by locally trivial fibrations
$\widetilde{\pi}_{b}~:~\widetilde{\mathcal {R}}~\to~\widetilde{W}$ and
$\widehat{\pi}_{b}~:~\widehat{\mathcal {R}}~\to~\widehat{W}$. Hence
$\mathfrak{g}_{0}~(\widetilde{\mathcal {R}}, \widetilde{\mathcal {F}})~=~0$,
$\mathfrak{g}_{0}(\widehat{\mathcal {R}},\widehat{\mathcal {F}})=0$, and
$\widetilde{W}=\widetilde{\mathcal {R}}/ \widetilde{\mathcal {F}}$,
$\widehat{W}=\widehat{\mathcal {R}}/\widehat{\mathcal {F}}$ are manifolds.

Since the fibrations $\widetilde{r}:\widetilde{M}\to B$ and $r:\widehat{M}\to B$
have the same base $B$,  each leaf of the foliation $(\widetilde{M}, \widetilde{F})$
is invariant respectively the group $\widehat{G}$, i.e.
$\widehat{G}\subset A_{L}(\widetilde{M}, \widetilde{F}).$ Therefore
$\widehat{\Gamma}\subset A_{L}(\widetilde{\mathcal {R}},\widetilde{\mathcal {F}})$
and the leaf spaces $\widetilde{\mathcal {R}}/ \widetilde{\mathcal {F}}=\widetilde{W}$
and $\widehat{\mathcal {R}}/\widehat{\mathcal {F}}= \widehat{W}$ are coincided, i.e.
$\widetilde{W}=\widehat{W}.$ Consequently, basic automorphism groups
$A_{B}(\widetilde{\mathcal {R}}, \widetilde{\mathcal {F}})$ and
$A_{B}(\widehat{\mathcal {R}},\widehat{\mathcal {F}})$ may be identified.
Further we put $A_{B}(\widetilde{\mathcal {R}}, \widetilde{\mathcal {F}})=
A_{B}(\widehat{\mathcal {R}},\widehat{\mathcal {F}}).$

According to the conditions of Theorem \ref{Th5}, $\Psi$ is a
discrete subgroup of the Lie group $Aut(B,\eta)$. Let $N(\Psi)$ be the
normalizer of $\Psi$ in the Lie group $Aut(B,\eta)\cong A^H(W, \widetilde{\beta})$.
Hence, $N(\Psi)$ is a closed Lie subgroup of the Lie group $Aut(B,\eta)$
and the quotient group $N(\Psi)/ \Psi$ is also a Lie group.

Let $\pi:{\mathcal {R}}\to M$ be the projection of the foliated
bundle over $(M, F)$. Due to the discreteness of
the global holonomy group $\Psi $,  the lifted  foliation $({\mathcal {R}}, {\mathcal {F}})$ is formed by
fibres of some locally trivial fibration $\pi_{b}:{\mathcal {R}}\to W$, which is called  the basic fibration.

Observe that there exists a map $\tau:\widehat{W}\rightarrow W$ satisfying
the equality $\tau\circ\widehat{\pi}_{b}=\theta \circ\pi_{b}$.
It is easy to show that $\tau:\widehat{W}\to W$ is a regular covering map with
the deck transformations group  ${\Phi}$, ${\Phi}\subset A^{H}(\widehat{W},\widehat{\beta})$,
which is naturally isomorphic to each of the groups $\Psi$, $G$ and $\Gamma$.

Denote by $\eta=(P(B,H),\omega)$  the Cartan geometry with the
projection $p:P\to B$ onto $B$ determined in the proof of
Theorem~\ref{Th2}. Remark that $\widehat{W}={P}$ is the space of the
$H$-bundle of the Cartan geometry $\eta$.

Since  $\kappa:\widehat{M}\to M$, $\theta:\widehat{\mathcal {R}}\to {\mathcal {R}}$ and $\pi: {\mathcal {R}}\to M$ are morphisms of
the following foliations
$\kappa:(\widehat{M},\widehat{F})\to (M,F)$,  $\theta:(\widehat{\mathcal {R}},\widehat{\mathcal {F}})\to ({\mathcal {R}},{\mathcal {F}})$ и $\pi:({\mathcal {R}},{\mathcal {F}}) \to (M,F)$ in the category of the foliations $\mathfrak F\mathfrak o\mathfrak l$, then maps $\widehat{\tau}:B\to M/F$
and $s:W\to W/H \cong M/F$ are defined, and the following diagram

\begin{equation}
	\label{D}
	\begin{CD}
		\xymatrix{
			P=\widehat{W}\ar@{->}[ddd]_{p}\ar@{->}[rrrr]^{\tau}&&&& W\ar@{->}[ddd]^{s}\\
			& \widetilde{\kappa}^{*}{\mathcal R}=\widetilde{\mathcal R}\ar@{->}[ul]^{\widetilde{\pi}_{B}}\ar@{->}[d]_{\widetilde{\pi}}\ar@{->}[r]_{\widehat{\theta}}&\kappa^{*}{\mathcal{R}}=\widehat{\mathcal{R}} \ar@{->}[ull]_{\widehat{\pi}_{B}}\ar@{->}[r]_{\theta}\ar@{->}[d]_{\widehat{\pi}}&{\mathcal {R}}\ar@{->}[ur]^{{\pi}_{B}}\ar@{->}[d]^{\pi}& \\
			& \widetilde{ M}\ar@{->}[ld]_{\widetilde{r}}\ar@{->}[r]^{\widehat{\kappa}} &\widehat{M}\ar@{->}[dll]_{r}\ar@{->}[r]^{\kappa}&M\ar@{->}[dr]^{q}&\\
			\widehat{M}/\widehat{F}= B\ar@{->}[rrrr]^{\widehat{\tau}}&&&& M/F  }\nonumber
	\end{CD}
\end{equation}
is commutative.

Due to Proposition \cite[Thm. 1]{SZ} there are the Lie group isomorphisms
$$\varepsilon: A_{B}(M, F)\to im(\varepsilon)\subset A^{H}(W,\widetilde{\beta})\,\, and\,\,$$
$$\widehat{\varepsilon}:A_{B}(\widetilde{M},\widetilde{F})=A_{B}(\widehat{M},
\widehat{F})\to im(\widehat{\varepsilon})\subset A^{H}(\widehat{W}, \widehat{\beta}).$$

Let us define a map $\Theta: im(\varepsilon)\to N(\Phi)/\Phi$ by the
following a way. Take any  $h\in im(\varepsilon)\subset
A^{H}(W,\widetilde{\beta})$. Denote the element
$\varepsilon^{-1}(h)\in A_{B}(M, F) $ by $f\cdot A_{L}(M, F)\in
A_{B}(M, F)$, where $f\in A(M, F)$. Since $\widetilde{\kappa}:\widetilde{M}\to M$
is the universal covering map there exists $\widetilde{f}\in Diff(\widetilde{M})$ 
lying over $f$ relatively $\widetilde{\kappa}$. It not difficult to see that
$\widetilde{f}\in A(\widetilde{M},\widetilde{F})$. Hence $\widetilde{f}\cdot
A_{L}(\widetilde{M}, \widetilde{F})\in A_{B}(\widetilde{M},
\widetilde{F})$. Consider $\widehat{h}:=\widehat{\varepsilon}(\widetilde{f}\cdot
A_{L}(\widetilde{M}, \widetilde{F}))\in im(\widehat{\varepsilon})\subset
A^{H}(\widehat{W},\widehat{\beta})$. The direct check shows that
$\widehat{h}$ lies over $h$ respectively $\tau$. Remind that $\Phi$ is
the deck transformation group  of the covering map
$\tau:\widehat{W}\to W$. Applying the 
Proposition~$\ref{pr8.2}$, we get that $\widehat{h}\in N(\Phi)$, hence
the set of all automorphisms in $im(\widehat{\varepsilon})$ lying
over $h$ is equal to the set of transformations from the class
$\widehat{h}\cdot\Phi.$ Let us put $\Theta(h):=
\widehat{h}\cdot\Phi\in N(\Phi)/\Phi.$ It is easy to check that the
map $\Theta: im(\varepsilon)\to N(\Phi)/\Phi$ is a group
monomorphism.

The effectiveness of the Cartan geometry $\eta=(P(B,H),\omega)$ on
$B,$ where $P=\widehat{W},$ implies the existence of the Lie group
isomorphism $\sigma:A^{H}(\widehat{W},\widehat{\beta})\rightarrow
Aut(B,\eta)$. Observe that
$\sigma({\Phi})=\Psi$ and $\sigma( N(\widetilde{\Phi}))=N(\Psi)$,
hence there exists the inducted Lie group isomorphism
$\widetilde{\sigma}: N({\Phi})/{\Phi}\to N(\Psi)/\Psi$. Thus, the
composition of the Lie group monomorphisms
\begin{eqnarray}
	\delta:=\widetilde{\sigma}\circ {\Theta}\circ\varepsilon:A_{B}(M,F)\rightarrow N(\Psi)/\Psi\nonumber
\end{eqnarray}
is the required Lie group monomorphism. Due to uniqueness  of the Lie group
structure in $A_{B}(M, F)$, in conforming with \cite[Thm.~1]{SZ},
the image $im(\delta)$ is an open-closed subgroup of the Lie group $N({\Psi})/{\Psi}$.
\qed
\section{Basic automorphism groups of  Cartan foliations\\ with an integrable Ehresmann connection}
\subsection{Proof of Theorem \ref{Th7}} 1. 
According to the conditions of the theorem  being proved, $(M,F)$ has  an intagrable Ehresmann connection  $\mathfrak M$.
In this case, distribution $\mathfrak M$ is integrable. In this case, there is $q$-dimensional
foliation $(M,F^t)$ such that $TF^t=\mathfrak M.$ 

Let $\widetilde{\kappa}: \widetilde{M}\to M$
be the universal covering map. According to the decomposition theorem belonging to
S.~Ka\-shi\-wabara \cite{Kas}, the universal covering manifold  $\widetilde{M}$ is equal to the product of manifolds $\widetilde{M} = \widetilde{Q}\times B$, where  $\widetilde{Q}$ is the universal covering manifold for any leaf of the foliation $(M, F)$, and $B$ is the universal covering manifold for any leaf  of the foliation $(M,F^t)$. The induced foliations $\widetilde{F}=\widetilde{\kappa}^*F=\{\widetilde{Q}\times\{y\}\,|\,
y\in B\}$, $\widetilde{F}^t=\widetilde{\kappa}^*F^t=\{\{z\}\times
B\,|\, z\in \widetilde{Q}\}$ are defined. Therefore, (M, F) is covered by fibration $\widetilde{s}: \widetilde{Q}\times B\to B$. In this case, by the same way as in the proof of
Theorem \ref{Th2}, the Cartan geometry $\eta$ is induced on $B$ such
that $(M,F)$ becomes an $(Aut(B,\eta),B)$-foliation.

2. Let $\Psi$ be the global holonomy group of this foliation.
Suppose now that the normalizer $N(\Psi)$ is equal to the
centralizer $Z(\Psi)$ of the group $\Psi$ in the group $Aut(B,\eta)$.

Let us fix points $x_0\in M$ and $(z_0,y_0)\in\widetilde{\kappa}^{-1}(x_0)\in\widetilde{M}$.
Then the fundamental group $\pi_1(M,x_0)$ acts on the universal covering space
$\widetilde{M} = \widetilde{Q}\times B$ as the deck transformation group
$\widetilde{G}\cong\pi_1(M,x_0)$ of $\widetilde{\kappa}$. Since $\widetilde{G}$
preserves both the inducted foliations $(\widetilde{M},\widetilde{F})$ and
$(\widetilde{M},\widetilde{F}^t)$, then every
$\widetilde{g}\in\widetilde{G}$ may be written in the form $\widetilde{g}=(\psi^t,\psi)$,
where $\psi^t$ generates a subgroup $\Psi^t$ in $Diff(\widetilde{Q})$,
$\psi\in\Psi$, and $\widetilde{g}(z,y)=(\psi^t(z),\psi(y))$,
$(z,y)\in\widetilde{Q}\times B.$  The maps $\widetilde{\chi}:
\widetilde{G}\to\Psi: \widetilde{g}\to\psi$ and $\widetilde{\chi}^t: \widetilde{G}\to\Psi^t:
\widetilde{g}\to\psi^t$ are the group epimorphisms.

Let $h$ be any element from
$N(\Psi)/\Psi$. Since $N(\Psi) = Z(\Psi)$, we have the following chain of equalities
$$\widetilde{g}\circ (id_{\widetilde{Q}}, h) =({\psi}^t,\psi)\circ
(id_{\widetilde{Q}}, h)=({\psi}^{t} \circ id_{\widetilde{Q}}, \psi \circ h)=\\$$
$$(id_{\widetilde{Q}} \circ {\psi}^{t}, h\circ \psi )= (id_{\widetilde{Q}}, h)\circ
({\psi}^t,\psi) = (id_{\widetilde{Q}}, h)\circ \widetilde{g}$$
for any $\widetilde{g}= ({\psi}^t,\psi)\in\widetilde{G}$, i.e.
$\widetilde{G}\cdot(id_{\widetilde{Q}},h)=(id_{\widetilde{Q}},h)\cdot\widetilde{G}$.
Therefore, by the Proposition~\ref{pr8.2}, for the deck transformation
group $\widetilde{G}$, there exists $\widetilde{h}\in Diff(M)$
such that $(id_{\widetilde{Q}},h)$ lies over $\widetilde h$ respectively to
$\widetilde{\kappa}:\widetilde{M}\to M$.

Taking into account that $(id_{\widehat{Q}},h)\in A(\widetilde{M},\widetilde{F})$, it is not difficult
to check that $\widetilde{h}~\in~A(M,F)$. Hence, $\varepsilon(\widetilde{h}\cdot A_{L}(M,F) )= h.$
This means that $\varepsilon: A_B(M,F)\to N(\Psi)/\Psi$ is surjective.
Thus, $\varepsilon$ is a the group isomorphism.

Since $N(\Psi)$ is a closed Lie subgroup of the automorphism Lie group $Aut(B,\eta)$, and $\Psi$ is a discrete subgroup of $N(\Psi)$, the quotient group $N(\Psi)/\Psi$ is a Lie group. Therefore, the group isomorphism $\varepsilon$ induces a Lie group structure in $A_B(M,F)$ such that $\varepsilon: A_B(M,F)\to N(\Psi)/\Psi$ becomes a Lie group isomorphism. According to \cite[Thm.~1]{SZ}, the Lie group structure in $A_B(M,F)$ is unique. \qed

\subsection{Example of finding a basic automorphism group}\label{Ex} Let $\mathbb{S}^{q}$ be a $q$-dimensional standard sphere, where $q\geq 3$. We identify $\mathbb{S}^{q}$  with $\mathbb{R}^q\cup\{\infty\}$, where $\{\infty\}$ is the point  at infinity. Define the transformation $\psi:\mathbb{S}^{q}\cong \mathbb{R}^q\cup\{\infty\} \to \mathbb{S}^{q}$ by equality $\psi(z)=\lambda z $  for any $ z\in \mathbb{S}^{q}\cong \mathbb{R}^q\cup\{\infty\}$, where $\lambda$ is a real number, and $0<\lambda<1$. We denote by $Conf(\mathbb{S}^{q})$   the Lie group of all conformal transformations of the sphere $\mathbb{S}^{q}$.

 Let $\Psi=<\psi>$ be the subgroup of the group $Conf(\mathbb{S}^{q})$ generated by $\psi$, and $\Psi$ is isomorphic to the group of integers $\mathbb{Z}$. Define the action of the group
$\mathbb{Z}$ on the product of manifolds $\mathbb{R}^{1}\times
\mathbb{S}^{q}$ by the equality $n(t,z)=(t-n,\psi^{n}(z))$ for any $n\in \mathbb{Z}$, $(t,z)~\in~\mathbb{R}^{1}\times
\mathbb{Z}$. This action is free and
 properly discontinuous. Therefore,  the mani\-fold of orbits $M~=~\mathbb{R}^1\times_{\mathbb{Z}}\mathbb{S}^{q}$ is defined. Denote by   $f:\mathbb{R}^{1}\times \mathbb{S}^{q} \to M$ the quotient map. Fix a point $(t_0,z_0)\in \mathbb{R}^{1}\times \mathbb{S}^{q}$, put $x_{0}=f(t_0,z_0)\in M$. Then the fundamental group $\pi_1(M,x_0)$ acts on the universal covering space
$\mathbb{R}^{1}\times \mathbb{S}^{q}$ as the deck transformation group
$\widetilde{G}\cong\pi_1(M,x_0)$ of $f$. Since the action  $\widetilde{G}$ preserves the structure of the product $\mathbb{R}^{1}\times \mathbb{S}^{q}$, then two foliations 
$(M,F)$ and $(M,F^t)$, covered by trivial fibrations $pr_{2}: \mathbb{R}^{1}\times \mathbb{S}^{q} \to \mathbb{S}^{q}$ and  $pr_{1}: \mathbb{R}^{1}\times \mathbb{S}^{q} \to \mathbb{R}^{1}$ respectively,  are defined.
Let us denote by
$\chi:\mathbb{R}^{1}\to \mathbb{S}^{1}=\mathbb{R}^{1}/\mathbb{Z}$  and $\nu: \mathbb{S}^{q}\to \mathbb{S}^{q}/\Psi$ the quotient maps onto the orbit spaces. Let $r:M\to M/F$ be the quotient map onto the leaf space. Observations
show that the topological spaces $M/F$ and $\mathbb{S}^{q}/\Psi$ are homeomorphic
and satisfy the commutative diagram

	\begin{equation}
	\label{D}
	\begin{CD}
		\xymatrix{
			\mathbb{R}^{1}\ar@{->}[d]_{\chi}\ar@{<-}[r]^{pr_{1}}&
			\mathbb{R}^{1}\times \mathbb{S}^{q}\ar@{->}[r]^{\,\,\,\,\,\,\,\,pr_{2}}\ar@{->}[d]_{f}&\mathbb{S}^{q}\ar@{->}[d]_{\nu}&\\
			\mathbb{S}^{1}\ar@{<-}[r]_{\,\,\,\,\,\,p}& M\ar@{->}[r]_{r}& \,\,\,\,\mathbb{S}^{q}/\Psi\cong M/F,&
		}\nonumber
	\end{CD}
\end{equation}  
where $p:M\to \mathbb{S}^{1}$ is the projection of the locally trivial fibration transforming the orbit $\mathbb{Z}.(t,z)$ of a point  $(t,z)\in \mathbb{R}^{1}\times\mathbb{S}^{q}$,   considered as a point from $M$, into the orbit  $\mathbb{Z}.{t}$ of a point  $t\in \mathbb{R}^{1}$, considered as a  point of the circle $\mathbb{S}^{1}$.  Since the manifold $M$ is the space of a locally trivial fibration $p:M\to \mathbb{S}^{1}$  over the circle $\mathbb{S}^{1}$ with a compact standard fiber $\mathbb{S}^{q}$, then $M$ is compact.

The distribution $\mathfrak M$ tangent to $(M, F^t)$, is an integrable Ehresmann connection for the foliation $(M, F)$. The foliation $(M, F)$ has two compact leaves $L_1$ and $L_2$ which are  diffeomorphic to the circle $\mathbb{S}^{1}$. Every other leaf $L$ of $(M, F)$ is diffeomorphic to $\mathbb{R}^1$, and its closure $\overline{L}$ is equal to the union $L\cup L_{1}\cup L_{2}$.
We emphasize that $(M, F)$ is a proper conformal foliation, which can be regarded as a Cartan foliation of type $(G, H)$, where $G=Conf(\mathbb{S}^{q})$ and $H$ is a stationary subgroup of the group $Conf(\mathbb{S}^{q})$ at some point in $\mathbb{S}^{q}$.

As is known, $H\cong CO(q)\ltimes \mathbb{R}^{q}$  is a semidirect product of a conformal group  $CO(q)\cong \mathbb{R}^{+}\times O(q)$ and a normal abelian subgroup $\mathbb{R}^{q}$. Note that $\Psi$ is the global holonomy group of the foliation $(M, F)$, and $\Psi$ is a discrete subgroup of the Lie group $Conf(\mathbb{S}^{q})$.

The direct check shows that the normalizer of the group $\Psi$ is equal to $N(\Psi)=\mathbb{R}^{+}\times O(q)$, and $N(\Psi)$ 
coincides with the centralizer $ Z(\Psi) $. Applying Theorem~\ref{Th7}, we obtain
that the group of basic automorphisms $ A_{B}(M,F)$ is a Lie group isomorphic to the quotient group $ N(\Psi)/\Psi\cong U(1)\times O(q)$, where $U(1)\cong \mathbb{S}^{1}$. Thus, the Lie group of basic automorphisms $ A_{B}(M,F) $ is isomorphic to the product of Lie groups $U(1)\times O(q) $.

\paragraph{Acknowledgements} The author would like to thank  N.I. Zhukova for helpful discus\-sions and comments.
The  work was supported by Laboratory of Dynamical Systems and
Applications NRU HSE, of the Ministry of science and higher education of the RF grant ag. no. 075-15-2019-1931

\end{document}